\magnification = 1200

\def\1b{{\bf 1}}

\def\eb{{\bf e}}

\def\gb{{\bf g}}

\def\mb{{\bf m}}
\def\nb{{\bf n}}

\def\Ab{{\bf A}}

\def\Cb{{\bf C}}

\def\Nb{{\bf N}}

\def\Hb{{\bf H}}
\def\Pb{{\bf P}}

\def\Rb{{\bf R}}

\def\Ub{{\bf U}}

\def\Aut{\rm{Aut}\,}

\def\Tor{{\rm{Tor}\,}}
\def\Hom{{\rm{Hom}\,}}
\def\Homc{{\cal H}{\it om}}
\def\Cohc{{\cal C}{\it oh}}
\def\End{{\rm{End}\,}}

\def\Ac{{\cal A}}
\def\Bc{{\cal B}}

\def\Ec{{\cal E}}
\def\Fc{{\cal F}}
\def\Gc{{\cal G}}
\def\Hc{{\cal H}}

\def\Kc{{\cal K}}
\def\Lc{{\cal L}}
\def\Mc{{\cal M}}
\def\Nc{{\cal N}}
\def\Oc{{\cal O}}

\def\Rc{{\cal R}}
\def\Sc{{\cal S}}
\def\Tc{{\cal T}}

\def\Wc{{\cal W}}

\def\qed{\hfill $\sqcap \hskip-6.5pt \sqcup$}        % White box               

\overfullrule=0pt                                    % No black boxes

\baselineskip =13pt

\baselineskip =13pt
\overfullrule =0pt

\centerline{\bf KLEINIAN SINGULARITIES, DERIVED CATEGORIES}
\centerline{\bf AND HALL ALGEBRAS}

\vskip2mm

\centerline{\bf M. Kapranov and E. Vasserot}

\vskip2cm

The Kleinian singularities ${\bf C}^2/G$ associated to finite
subgroups $G\i SL_2({\bf C})$, are of fundamental importance
in algebraic geometry, singularity theory and other 
branches of mathematics. Despite the very classical nature of the subject,
new remarkable properties continue to be discovered. One such discovery
was the McKay correspondence [9] and its interpretation by Gonzalez-Springberg
and Verdier [3] in terms of the minimal resolution ${{\bf C}^2//G}$. Their 
results give identifications
$$K_0(
{{\bf C}^2//G}) \quad \simeq \quad {\rm Rep}(G) \quad \simeq 
\quad \widehat{\bf h}_{\bf Z},$$
where $K_0$ is the Grothendieck group, Rep is the representation ring and
$ \widehat{\bf h}_{\bf Z}$ is the root lattice of the affine Lie algebra
(of type A-D-E) associated to $G$. 

\vskip .1cm

Our first goal in this paper is to extend the above results
by describing the derived category of coherent sheaves on
${{\bf C}^2//G}$, instead of just $K_0$. Theorem 1.4
identifies it with the derived category of $G$-equivariant
${\bf C}[x,y]$-modules, i.e., of modules over the
crossed product algebra ${\bf C}[x,y] [G]$. It is surprising that
such a basic fact has not been noticed before. Our approach can be seen
as a refinement, in a purely algebraic setting, of the techniques of
Kronheimer and Nakajima [7] in that we get rid of Dolbeault
complexes with growth conditions at infinity, stability
conditions for vector bundles and so on. 

\vskip .1cm

We then define {\bf H}, an Euler-characteristic version of the Hall algebra
 [12]
of the category of coherent sheaves on ${{\bf C}^2//G}$ and apply the 
constructed equivalence to exhibit
a subalgebra in {\bf H} isomorphic to $U({\bf g}^+_G)$.
Here ${\bf g}^+_G$ is the nilpotent part of the finite-dimensional
Lie algebra (of type A-D-E) corresponding to $G$. As a consequence,
we get a result about any algebraic surface $S$ equipped with
a configuration $C=\bigcup_i {\bf P}^1_i$ of (-2)-curves intersecting
transversally. Namely, taking the intersection graph of the ${\bf P}^1_i$
as a Dynkin graph, we get a possibly infinite-dimensional Kac-Moody
Lie algebra, and the theorem is that the positive part of this algebra
acts in the space of functions on isomorphism classes of coherent
sheaves in $S$. This partly extends the results of Nakajima [10]
to a wider geometric context. 

\vskip .1cm

The work of the first author was partly supported by an NSF grant and a
large part of it was carried out during a visit to Universit\'e
Cergy-Pontoise, whose hospitality is gratefully acknowledged. 
The work of the second author was partly supported by EEC grant 
no. ERB FMRX-CT97-0100.

\vfill\eject

\centerline {\bf \S 1. Equivalence of derived categories.}

\vskip3mm 

\noindent {\bf 1.1.}
For a smooth algebraic variety $X$ and an integer $m\geq 0$ we denote
by $X^{(m)}$ the $m$th symmetric power of $X$ and by
$X^{[m]}$ the Hilbert scheme parametrizing 0-dimensional
subschemes $\xi\i X$ of length $m$. Given such a $\xi$, 
the corresponding point of $X^{[m]}$ is denoted $[\xi]$. 
We denote by $\cal E$ the tautological $m$-dimensional bundle on $X^{[m]}$,
whose fiber at $[\xi]$ is $H^0(\xi, {\cal O})$. 
Let $G$ be a finite group acting on $X$.
Let $m=|G|$. The quotient $X/G$ can be viewed
as a closed subvariety in $X^{(m)}$.
Suppose that the action is free on an open $G$-invariant set $U\i X$.
Define the Hilbert quotient $X//G$ as the closure of
$U/G$ in $X^{[m]}$, cf. [6]. The Chow morphism
$X^{[m]}\to X^{(m)}$ gives a map $p\,:\, X//G \to X/G$.
Let $\Sigma\i (X//G)\times X$ be the incidence subscheme.
Let $p_1: (X//G)\times X\to X//G$, $p_2: (X//G)\times X\to X$,
be the projections and $q_1, q_2$ be the restrictions of $p_1, p_2$
to $\Sigma$. The restriction of the tautological bundle to $X//G$ 
is denoted again by $\cal E$; it is a bundle of $G$-modules 
isomorphic to the regular representation. 
If $\rho,V$ are representations of $G$, with $\rho$ irreducible,
we set $V_\rho=\Hom_G(\rho,V)$.

\vskip .3cm

\noindent {\bf 1.2.} 
Let $\Cohc_G(X)$ be the category of $G$-equivariant coherent sheaves 
on $X$, and $\Cohc(X//G)$ be the category of coherent sheaves on $X//G$. 
Define two functors $\Phi: D^b(\Cohc_G(X))\to D^b(\Cohc(X//G))$  and
$\Psi: D^b(\Cohc(X//G))\to D^b(\Cohc_G(X))$ by
$$\Phi({\cal F})=(Rq_{1*}Lq_2^*{\cal F})^G = 
(Rp_{1*}(p_2^*{\cal F}\otimes^L {\cal O}_\Sigma))^G
\quad{\rm and}\quad
\Psi(\Gc) = Rp_{2*}R\Homc(\Oc_\Sigma, p_1^*\Gc).$$
The functors $\Psi$ and $\Phi$ are adjoint, i.e.
$\Homc(\Phi({\cal F}), {\cal G}) \simeq \Homc({\cal F}, \Psi({\cal G})).$

\vskip .3cm

\noindent {\bf 1.3.} 
Let $Z,W$ be algebraic varieties,
equipped with actions of finite groups $G$, $H$ respectively. 
Let $p_W, p_Z$ be the projections from $Z\times W$ to $W,Z$ respectively.
Let $\cal L$ be an object of $D^b(\Cohc_{G\times H}(Z\times W))$, 
such that each of the cohomology sheaves
${\cal H}^i({\cal L})$ has proper support with respect to $p_W$. 
Taking  $\cal L$  as a ``kernel" defines a functor
$$F_{\cal L}: D^b(\Cohc_G(Z))\to D^b(\Cohc_H(W)), \quad {\cal F}\to
(Rp_{W*}(p_Z^*{\cal F}\otimes^L {\cal L}))^G.$$
If $Z,W,T$ are varieties with actions of groups
$G,H,K$ respectively, and $\Lc\in D^b (\Cohc_{G\times H}(Z\times W))$,
$\Mc\in D^b (\Cohc_{H\times K}(W\times T))$,
the composition $F_{\cal M}\circ F_{\cal L}$ is isomorphic to
$F_{{\cal M}*{\cal L}}$ where
$${\cal M}*{\cal L}=(Rp_{13*}(p_{12}^*{\cal L}\otimes^L p_{23}^*{\cal M}))^H,$$
and the $p_{ij}$ are the projections of $Z\times W\times T$ onto the pairwise 
products. In the case (1.2) of the product $(X//G)\times X\times (X//G)$ let 
$\Sigma_{12}=p_{12}^{-1}(\Sigma)$ and $\Sigma_{23}^t = p_{23}^{-1}(\Sigma^t),$
where $\Sigma^t\i X\times (X//G)$ is the transpose variety of $\Sigma$.

\vskip3mm

\noindent{\bf Proposition.} {\it In the situation of (1.2) the
 following is true:}
\hfill\break
{\it (a) The composition $\Phi\Psi$ has the kernel 
$${\cal L}=(Rp_{13*}(R\Homc({\cal O}_{\Sigma_{12}},  
{\cal O}_{\Sigma_{23}^t})))^G\in D^b(\Cohc(X//G\times X//G)).$$
(b) The composition $\Psi\Phi$ has the kernel 
$${\cal M} = Rs_{13*}(R\Homc({\cal O}_{\Sigma_{12}^t},
{\cal O}_{\Sigma_{23}}))\in D^b(\Cohc_{G\times G}(X\times X)),$$
where $s_{ij}$ are the projections of $X\times (X//G)\times X$ to pairwise 
products. 
}\qed

\vskip3mm

\noindent{\bf Remark.}
Note that that the  restriction of the projection
$$p_{13}: \Sigma_{12}\cap\Sigma_{23}^t\to (X//G)\times (X//G)$$ 
is a finite morphism, so we can replace $Rp_{13*}$ by $p_{13*}$. 

\vskip3mm

\noindent
Thus, the functors $\Phi, \Psi$ are mutually inverse equivalences of 
categories
if and only if \hfill\break
(a) the kernel $\cal L$ is quasi-isomorphic to the structure sheaf of
the diagonal on $(X//G)\times (X//G)$,\hfill\break
(b) the kernel $\cal M$ is quasi-isomorphic to the structure sheaf of diagonal 
$\Delta\i X\times X$ tensored with the regular representation 
$\Rb=\Cb[G]$ of $G$. 

\vskip3mm

\noindent {\bf 1.4.} 
Suppose now that $X$ is a smooth surface. Then $X^{[m]}$ is smooth and
 $X//G$ is an irreducible component of
the fixed point set of the $G$-action on $X^{[m]}$, so it is also smooth.
Thus $p$ is a resolution of singularities of $X/G$.
The following theorem is the main result of the section.

\proclaim Theorem. Let $X$ be a surface equipped with a holomorphic 
symplectic form $\omega$, and suppose that the $G$-action on $X$ preserves 
$\omega$. Then $\Phi$ and $\Psi$ are mutually inverse  equivalences of 
categories.

\vskip3mm

\noindent  From the definitions it is clear that all the cohomology 
sheaves of $\cal L$ are supported on the $([\xi],[\eta])$ such that
${\rm supp}(\xi)\cap{\rm supp}(\eta)\neq\emptyset,$
and the cohomology sheaves of $\cal M$ are supported on
$\bigcup_{g\in G} (1\times g)(\Delta).$
Thus it is enough to work in the neighborhoods 
of fixed points of subgroups of 
$G$, where the action can be replaced by the linear one. 
We now assume this to be the case. 

\vskip3mm

\noindent {\bf 1.5.} 
Given a finite subgroup $G\subset SL_2(\Cb)$ let $\tau$ be
the natural representation of $G$ on $\Cb^2$.
For any irreducible representations $\pi, \rho$ of $G$ let 
$m_{\pi\rho}$ be the multiplicity of $\pi$ in $\rho\otimes_\Cb\tau$.
We now assume $X=\tau$ and $G\i SL_2({\bf C})$. 
Set $\Ab={\bf C}[x,y]$. Let $M$ be any $\Ab$-module. 
Then, we have the (Koszul) free resolution of $M$ by
$$\Ab\otimes_\Cb M\buildrel (x,y)\over\to (\Ab\otimes_\Cb M)^{\oplus 2} 
\buildrel\pmatrix{y\cr-x}\over\to\Ab\otimes_\Cb M.$$
Constructing this resolution simultaneously for $M=\Oc_\xi$ and all
$[\xi]\in X//G$, we get a resolution of ${\cal O}_\Sigma$ 
by the complex on $(X//G)\times X$ :
$$\Kc=\{p_1^*\Ec\to
(p_1^*\Ec)^{\oplus 2}\to p_1^*\Ec\}.$$
Thus 
$$\matrix{
\Phi(\Fc)&=(Rp_{1*}(p_2^*\Fc\otimes^L\Kc))^G\hfill\cr 
&=\{ (Rp_{1*}p_2^*\Fc)\otimes\Ec\to 
(Rp_{1*}p_2^*\Fc)\otimes\Ec^{\oplus 2}
\to (Rp_{1*}p_2^*\Fc)\otimes\Ec\}^G\hfill.}$$
Note that $Rp_{1*}p_2^*\Fc$ is just the trivial bundle on $X//G$ with fiber
$\Gamma(X,\Fc)$ (the higher cohomology vanishes on $X=\tau$). 
Moreover, if $\Fc$ is a $G$-equivariant sheaf, then $\Gamma(X,{\cal F})$ 
is a $\Ab$-module with $G$-action and can be split
$\Gamma(X,\Fc) = \bigoplus_\pi \pi\otimes_\Cb\Gamma(X,{\cal F})_\pi,$
so $\Phi(\Fc)$ can be rewritten as
$$\Phi(\Fc)=\{\bigoplus_\pi\Gamma(X,\Fc)_\pi\otimes_\Cb\Ec_\pi\to 
\bigoplus_{\pi,\rho}\Gamma(X,\Fc)_\pi\otimes_\Cb\Ec_\rho^{m_{\pi\rho}}\to 
\bigoplus_\pi\Gamma(X,\Fc)_\pi\otimes_\Cb\Ec_\pi\}.$$
Similarly we get
$$\Psi(\Gc)=\{ 
R\Gamma(X//G,\Gc\otimes\Ec^*)\otimes_\Cb\Oc_\tau\to
R\Gamma(X//G,\Gc\otimes\Ec^*)\otimes_\Cb\Oc_\tau^2 \to
R\Gamma(X//G,\Gc\otimes\Ec^*)\otimes_\Cb\Oc_\tau
\}.$$

\noindent 
By using the Koszul resolution $\Kc$, 
we find that $\Lc$ is quasi-isomorphic to
$$
\Lc'=(Rp_{13*}(R\Homc(p_{12}^*\Kc,\Oc_{\Sigma_{23}^t})))^G$$
$$=\{\Ec^*\underline{\otimes}\Ec\to 
(\Ec^*\underline{\otimes}\Ec)^{\oplus 2}\to
\Ec^*\underline{\otimes}\Ec\}^G$$
$$=\{\bigoplus_\pi\Ec_\pi^*\underline{\otimes}\Ec_\pi\to
\bigoplus_{\pi,\rho}\Ec_\pi^*\underline{\otimes}
\Ec_\rho^{m_{\pi\rho}}\to
\bigoplus_\pi\Ec_\pi^*\underline{\otimes}\Ec_\pi\},$$
where $\underline{\otimes}$ is the external tensor product.
The quasi-isomorphism $\Lc'\simeq\Oc_\Delta$ is proved  in [11, Lemma 4.10].
By using the Koszul resolution for $\Oc_{\Sigma_{12}^t}$
we represent $\cal M$ by the quasi-isomorphic complex
$$\Mc'=\{\Oc_\tau\underline{\otimes}Rq_{2*}q_1^*{\cal E}^* 
\to\bigl(\Oc_\tau\underline{\otimes}
Rq_{2*}q_1^*{\cal E}^*\bigr)^{\oplus 2}\to
\Oc_\tau\underline{\otimes}Rq_{2*}q_1^*\Ec^*\}.$$
To show that $\Mc'$ is quasi-isomorphic to $\Oc_\Delta\otimes_\Cb\Rb$, 
it is enough to show that 
$Rq_{2*}(q_1^*\Ec^*)=\Oc_{\tau}\otimes_\Cb\Rb.$
Then $\Mc'$ will be identified with the tensor product of $\Rb$
and the Koszul resolution of the diagonal in $\Cb^2\times\Cb^2$:
$$\Oc\underline{\otimes}\Oc\to\Oc\underline{\otimes}\Oc^{\oplus 2}
\to\Oc\underline{\otimes}\Oc.$$
Applying the Koszul resolution one more time, 
we are reduced to the following fact.

\vskip3mm

\noindent{\bf Proposition.}
{\it We have
$\Gamma(X//G,\Ec^*\otimes\Ec)=\Rb[x,y]$ and
$H^i(X//G,\Ec^*\otimes\Ec)=0\quad\forall i>0.$
}

\vskip3mm

\noindent{\sl Proof:} 
Recall that $\Gamma(X//G,\Ec)=\Ab$ (see [3])
and thus for any $\Gc\in\Cohc(X//G)$ the space
$\Gamma(X//G,\Ec\otimes\Gc)$ is a $\Ab$-module. 
We define a morphism of $\Ab$-modules
$u:\Rb[x,y]\to\Gamma(X//G,\Ec^*\otimes\Ec)$ by
$$u(g) = g\in\Hom(\Ec,\Ec)=\Gamma(X//G,\Ec^*\otimes\Ec),
\qquad\forall g\in G.$$
Since the $G$-action on $\tau$ is free outside 0,
the space $\Gamma(X//G,\Ec^*\otimes\Ec)$
is a torsion free $\Ab$-module, and $u$
is an isomorphism outside $0\in\tau$.
Thus the first assertion follows from
the following lemma.

\vskip3mm

\noindent{\bf Lemma.}
{\it Let $u: M\to N$ be a homomorphism of 
$\Ab$-modules such that:\hfill\break
{(a)} $M$ is free, and $N$ has no torsion;\hfill\break
{(b)} $u$ is an isomorphism outside a point (in particular, it is injective).
\hfill\break 
Then $u$ is an isomorphism.
}

\vskip3mm

\noindent {\sl Proof:} Let us regard $M,N$ as coherent sheaves $\Mc,\Nc$ on 
$\Cb^2$,
let $U$ be the complement of the point, so that over $U$ the map $u:\Mc\to\Nc$ 
is an isomorphism. Let $j: U\to\Cb^2$ be the embedding. 
Then, taking the direct image, we find a homomorphism
$$j_*j^*u: j_*j^*\Mc\to j_*j^*\Nc.$$ 
It is an isomorphism since $j^*u$ is. Since $\Mc$ is free,
$j_*j^*\Mc=\Mc$. On the other hand, since $\Nc$ 
is torsion free, the natural map $\Nc\to j_*j^*\Nc$ 
is an embedding. But its image contains the image of
$j_*j^*u$, so it is an isomorphism. 
So $u=j_*j^*u$ is an isomorphism, as claimed. 
\qed

\vskip3mm

\noindent
To prove the second assertion of the proposition it is enough to show that
$R^ip_*(\Ec^*\otimes\Ec)=0$ if $i>0.$
The map $p$ is an isomorphism everywhere except for $0\in X/G$. 
The proposition follows from the vanishing of $H^i(p^{-1}(0),\Ec^*\otimes\Ec)$
if $i>0$, proved in Remark 2.1 below.
\qed

\vskip5mm

\centerline {\bf \S 2. Detailled analysis of the equivalence.}

\vskip3mm
 
\noindent
In this section we always assume that $X=\tau$
and $G$ is a finite subgroup in $SL_2(\Cb)$. 

\vskip3mm

\noindent {\bf 2.1.}
The quotient surface $X/G$ has an isolated singularity at $0$. 
The map $p$ is the minimal desingularization of $X/G$ (see [5]).
The vector bundles $\Ec_\pi$ on $X//G$ were 
introduced and studied in [3].  
Recall that finite subgroups $G\i SL_2(\Cb)$ are classified by
Dynkin diagrams of finite type A-D-E as follows.
Let $E=p^{-1}(0)$ and let $E_{red}$ be the reduced variety.

\vskip3mm

\noindent{\bf Proposition.}
{\it Let $\Gamma^0$ be the Dynkin diagram of $G$.
Let $\Gamma$ be the affine extension of $\Gamma^0$.\hfill\break
(a) Vertices of $\Gamma^0$ are in bijection with components of $E$. 
Two vertices are joined by an edge if and only if the corresponding
components intersect. This intersection is transverse and consists of one 
point.
\hfill\break
(b) Vertices of $\Gamma$ are in bijection with irreducible representations of 
$G$. The vertices corresponding to $\pi$ and $\rho$ are joined if and only if 
$m_{\pi,\rho}\neq 0$. 
}\qed

\vskip3mm

\noindent
Let $\Pb^1_\pi\i E_{red}$ be the component corresponding 
to a nontrivial irreducible representation $\pi$. 
Put $d_\pi={\rm {rk}}\,\Ec_\pi={\rm {dim}}\,\pi$.

\vskip3mm

\noindent{\bf Lemma.}
{\it The restriction of $\Ec_\pi$ to $\Pb^1_\rho$ is trivial
(isomorphic to $\Oc^{d_\pi}$)
if $\pi\neq\rho$ and is isomorphic to $\Oc(1)\oplus\Oc^{d_\pi - 1}$, if
$\pi=\rho$.} 

\vskip3mm

\noindent {\sl Proof:}
Consider the (infinite-dimensional) space $\Ab_\pi$. 
By construction of $\Ec_\pi$, we have a surjective map of sheaves on 
$\Pb^1_\rho$ :
$$\Ab_\pi\otimes_\Cb\Oc_{\Pb^1_\rho}\to\Ec_\pi|_{\Pb^1_\rho}.$$
This implies that in the splitting 
$\Ec_\pi|_{\Pb^1_\rho}\simeq\bigoplus_i\Oc(m_i)$
each $m_i\geq 0$. Since the degree is $\sum_i m_i$, 
our statement follows from the following result of [2] : 
the degree of the restriction of ${\cal E}_\pi$
to $\Pb^1_\rho$ is 0 if $\pi\neq\rho$ and 1 if $\pi=\rho$. 
\qed

\vskip3mm

\noindent{\bf Remark.}
By Lemma 2.1, we know that for all $\pi\neq\Cb$ there exist integers
$a,b,c,$ such that
$$(\Ec^*\otimes\Ec)|_{\Pb^1_\pi}=
{\cal O}(-1)^{\oplus a}\oplus {\cal O}(0)^{\oplus b}\oplus
{\cal O}(1)^{\oplus c}.$$
Hence $H^i(\Pb^1_\pi,\Ec^*\otimes\Ec)=0$ for all $i>0$.
The irreducible components $E_\pi$ of $E$ may be non reduced but have
the self intersection number $(-2)$. Therefore, for any vector
bundle $\cal F$ on $X//G$ the restriction ${\cal F}|_{E_\pi}
=
{\cal F}\otimes
{\cal O}_{E_\pi}$ has a filtration with quotients of the form
${\cal F}|_{{\bf P}^1_\pi}\otimes {\cal O}_{{\bf P}^1_\pi}(2j), j\geq 0$.
Applying this to 
${\cal F}= \Ec^*\otimes\Ec$, we find that its
restriction  to each component  $E_\pi$ has no higher cohomology.
As a consequence, we get
$$H^i(E,\Ec^*\otimes\Ec)=0,\qquad\forall i>0.$$

\vskip2mm

\noindent{\bf 2.2. Lemma.} {\it 
(a) If $[\xi]\in E$ then $\Tor_0(\Oc_\xi,\Oc_0)=\Cb,$
$\Tor_1(\Oc_\xi,\Oc_0)=\Cb\oplus\bigoplus_{[\xi]\in\Pb^1_\pi}\pi$ and
$\Tor_2(\Oc_\xi,\Oc_0)=\bigoplus_{[\xi]\in\Pb^1_\pi}\pi$.\hfill\break
(b) The line bundles on $\Pb^1_\pi$ formed by the
$\Tor_1(\Oc_\xi,\Oc_0)_\pi$ and the $\Tor_0(\Oc_\xi,\Oc_0)^G$
are isomorphic to $\Oc(-1)$ and $\Oc$ respectivelly.
}

\vskip3mm

\noindent {\it Proof:}
Let ${\bf m} = (x,y)\i\Ab$ and let $\nb\subseteq\mb$ 
be the ideal generated by $\mb^G$. 
If $[\xi]\in E$, then the ideal $I_\xi\i\Ab$ contains $\nb$
(the $G$-module $\Oc_\xi=\Ab/I_\xi$ is isomorphic to the
regular representation, in particular, $\dim \Oc_\xi^G=1$. 
If an invariant $f\in\nb$
does not lie in $I_\xi$, then $f \,\,{\rm mod}\,\, I_\xi \in\Oc_\xi$ gives
a $G$-invariant on $\Oc_\xi$ not proportional to 1, so the dimension of
$\Oc_\xi^G$ is at least 2, a contradiction).
If $W\i {\bf m}/{\bf n}$ is $G$-invariant set $I(W)=\Ab\cdot W+\nb.$

\proclaim Theorem [5]. 
Let $\pi$ be a nontrivial irreducible representation of $G$. Then :\hfill\break
(a) There exist two different irreducible submodules 
$\pi', \pi''\i {\bf m}/{\bf n}$, isomorphic to $\pi$ and such that if
$[\xi]\in\Pb^1_\pi-\bigcup_{\rho\neq \pi}\Pb^1_\rho$ then 
$I_\xi = I(W)$, where $W\i \pi'\oplus\pi''$ 
is a proper nonzero $G$-submodule. \hfill\break
(b) If $m_{\pi,\rho}\neq 0$ then the point $[\xi]\in\Pb^1_\pi\cap\Pb^1_\rho$
has the ideal $I_\xi=I(\pi'\oplus\rho'')$. 
\qed

\vskip3mm

\noindent 
Part $(a)$ for $\Tor_0$ and the second claim of part $(b)$ of Lemma
follow from
the equality $\Oc_\xi\otimes\Oc_0=\Cb$. Since
$$\Tor_i(\Oc_\xi,\Oc_0)=\Tor_{i-1}(I_\xi,\Oc_0), \quad i=1,2,$$
we get $\Tor_1(\Oc_\xi,\Oc_0)=I_\xi/\mb I_\xi.$
If $[\xi]\in E$ then $I_\xi=I(W)$ where $W\i\mb/\nb$ and
$W_\pi=\bigoplus_{[\xi]\in\Pb^1_\pi}\pi$.
Now observe that $W=I_\xi/(\mb I_\xi +\nb)$ and that
$((\mb I_\xi+\nb)/\mb I_\xi)_\pi=0$ for any $\pi\neq\Cb$ since 
$(\mb I_\xi+\nb)/\mb I_\xi\simeq\nb/(\nb\cap\mb I_\xi)$
is a quotient of $\nb/\mb\nb$.
Thus, if $\pi\neq\Cb$,
$\Tor_1(\Oc_\xi,\Oc_0)_\pi$ is 1-dimensional
for $[\xi]\in\Pb^1_\pi$ and vanishes otherwise.
It remains to study $\Tor_2(\Oc_\xi,\Oc_0)_\pi$ and $\Tor_1(\Oc_\xi,\Oc_0)^G.$ 
The Tor's in question are the cohomology of the Koszul complex
$$\Oc_\xi\to\Oc_\xi\otimes_\Cb\tau\to\Oc_\xi.$$
Since $\Tor_2(\Oc_\xi,\Oc_0)\subseteq\Ab/I_\xi$ 
and $(\Ab/I_\xi)^G=\Cb+I_\xi$, we get $\Tor_2(\Oc_\xi,\Oc_0)^G=0$.
Part $(a)$ follows since the equivariant Euler characteristic of the complex 
is 0 (note that $\Oc_\xi\otimes_\Cb\tau\simeq\Oc_\xi\otimes_\Cb\Cb^2$
since $\Oc_\xi$ is the regular representation of $G$). 
To see part $(b)$, notice that
$\Pb^1_\rho\simeq\Pb\Hom_G(\rho,\rho'\oplus\rho''),$
and that the bundle formed by the 
$\Tor_1(\Oc_\xi,\Oc_0)_\rho=(I_\xi/{\bf m}I_\xi)_\rho$
on this $\Pb^1$ is just $W\mapsto W_\rho$, i.e. $\Oc(-1)$. 
\qed

\vskip3mm

\noindent{\bf 2.3.}
Any finite-dimensional representation $V$ of $G$
gives rise to two equivariant sheaves on $X$ : 
the skyscraper sheaf $V^!$ whose fiber at 0 is $V$ and all the other fibers 
vanish, and the locally free sheaf $\tilde V=V\otimes_\Cb\Oc_X$. 
If $\pi$ is irreducible, then $\Phi(\tilde\pi)=\Ec_\pi$. 
The sheaf $\pi^!$ is quasi isomorphic to the (Koszul) complex 
$$\tilde\pi\otimes_\Cb\Lambda^2\tau\to
\tilde\pi\otimes_\Cb\tau\to\tilde\pi.$$

\vskip3mm

\noindent{\bf Theorem.} {\it 
We have $\Phi(\Cb^!)=\Oc_E$ and $\Phi(\pi^!)=\Oc_{\Pb^1_\pi}(-1)[1]$
if $\pi\neq\Cb$ .}

\vskip3mm

\noindent{\sl Proof:} Let $\Fc=\Phi(\pi^!)$. Let us view it as the complex
of locally free sheaves on $X//G$ (see (1.5))
$$\Fc=\{\Ec_\pi\to\bigoplus_\rho\Ec_\rho^{m_{\pi\rho}}\to\Ec_\pi\}.$$
Then for $[\xi]\in X//G$
$$\Tor_i(\Oc_\xi,\Oc_0)_\pi=H^{-i}([\xi],\Fc_{[\xi]}).\leqno(2.3.1)$$
Recall that we have a spectral sequence
$$\Tor_i(\Hc^{-j}(\Fc),\Oc_{[\xi]})\Rightarrow H^{-i-j}([\xi],\Fc_{[\xi]}).
\leqno(2.3.2)$$
Lemma 2.2, $(2.3.1)$ and $(2.3.2)$ imply that
$\Hc^0(\Fc)=0$ if $\pi\neq\Cb.$
Moreover, if $\pi=\Cb$, if $\Lambda$ is any ring and if 
$[\xi]\in (X//G)(\Lambda)$ is the $\Lambda$-point corresponding to the 
subscheme $\xi\subset X\times{\rm Spec}(\Lambda)$ then
$$\Hc^0(\Fc)_{[\xi]}={\rm Tor}_0(\Oc_\xi,\Oc_0)^G=
{\rm Tor}_0(\Oc_{\bar\xi},\Oc_{\bar 0})=\Oc_{E,[\xi]},$$
where $\bar\xi$ and $\bar 0$ are the projection of $\xi$ and $0$ in $X/G$.
Thus, $\Hc^0(\Fc)=\Oc_{E}.$ 

By (2.3.1), (2.3.2) and Lemma 2.2, if $\pi\neq\Cb$ then 
$\Hc^{-1}(\Fc)$ is supported on $\Pb^1_\pi$ and its restriction
onto $\Pb^1_\pi$ is 
$$\Hc^{-1}(\Fc)/I_{\Pb^1_\pi}\Hc^{-1}(\Fc)\quad\simeq\Oc_{\Pb^1_\pi}(-1).$$
This does not yet imply that $\Hc^{-1}(\Fc)$ is actually a sheaf
on $\Pb^1_\pi$ rather
than on some infinitesimal neighborhood. For this, we need to show that 
$\Hc^{-1}(\Fc)$ is annihilated by the sheaf of ideals 
$I_{\Pb^1_\pi}$. Observe that in the 
category of $G$-equivariant sheaves we have $\End_G(\pi^!)=\Cb$, so 
all the endomorphisms of $\Fc$ in $D^b(\Cohc(X//G))$
are scalar. But $\Pb^1_\pi$ is a $(-2)$-curve, so it possesses a lot of 
functions regular in an entire neighborhood of $\Pb^1_\pi$ and vanishing on 
$\Pb^1_\pi$. So if $\Hc^{-1}(\Fc)$ actually is not annihilated by 
$I_{\Pb^1_\pi}$, there will be a section $f$ of $I_{\Pb^1_\pi}$ on some 
neighborhood of $\Pb^1_\pi$ which is not annihilating $\Hc^{-1}(\Fc)$.
Since all the homology of $\Fc$ is supported on $\Pb^1_\pi$, the
multiplication by such an  $f$
defines an endomorphism of $\Fc$ in the derived category.
This endomorphism is not scalar, since the
induced endomorphism on $\Hc^{-1}(\Fc)$ is not scalar. This contradicts the
assumption. So we have established that $\Hc^{-1}(\Fc)=\Oc_{\Pb^1_\pi}(-1)$
if $\pi\neq\Cb$. If $\pi=\Cb$ then $\Tor_1(\Oc_\xi,\Oc_0)^G=\Cb$ and
$\Tor_1(\Hc^0(\Fc),\Oc_{[\xi]})=\Oc(-E)_{[\xi]}=\Cb$. Since
$\Tor_1(\Oc_\xi,\Oc_0)^G$ is composed from $\Tor_1(\Hc^0(\Fc),\Oc_{[\xi]})$
and $\Tor_0(\Hc^{-1}(\Fc),\Oc_{[\xi]})$ we have $\Hc^{-1}(\Fc)=0$.

Finally, $\Hc^{-2}(\Fc) = 0$. More precisely, if $\pi\neq\Cb$ then
$\Tor_2(\Oc_\xi,\Oc_0)_\pi=\Cb$ is composed from
$\Tor_1(\Hc^{-1}(\Fc),\Oc_{[\xi]})$ (which is nonzero) and 
$\Hc^{-2}(\Fc)\otimes\Oc_{[\xi]}$. 
If $\pi=\Cb$ then $\Tor_2(\Oc_\xi,\Oc_0)^G=0$. 
\qed

\vskip1cm

\centerline {\bf \S 3. Hall algebras and double quivers.}

\vskip 5mm

\noindent {\bf 3.1.} Let us describe a version of the Hall algebra
construction [12] based on Euler characteristic. 
Let $\Ac$ be a $\Cb$-linear Abelian category of finite type
(i.e. the extension groups between pairs of objects in $\Ac$
are finite dimensional). If $A,B,C$ are three objects of $\Ac$, the set 
$G_{AB}^C=\{A'\subseteq C:\,\, A'\simeq A,\,C/A'\simeq B\}$
has the structure of a complex variety. 
To see this, let ${\rm Com}_{AB}^C$ and $E_{AB}^C$ be the set of all complexes
and all exact sequences respectively of the form
$0\to A\buildrel \alpha\over \to C\buildrel\beta\over\to B\to 0.$
Clearly, ${\rm Com}_{AB}^C$ is a closed algebraic subvariety in the affine 
space
$\Hom(A,C)\oplus\Hom(C,B)$,
and $E_{AB}^C$ is a Zariski open subset in ${\rm Com}_{AB}^C$. 
Now, the algebraic group $\Aut(A)\times\Aut(C)$ acts on $E_{AB}^C$ freely.
The quotient is therefore equipped with a structure of a complex variety. 
But as a set, this quotient is nothing but $G_{AB}^C$.
Since the heart of a triangulated category is stable by
extensions, we get

\vskip3mm

\noindent{\bf Proposition.}
{\it Let $\Ac, \Bc$ be two Abelian categories as above, and 
$F: D^b(\Ac)\to D^b(\Bc)$ be an equivalence of triangulated categories. 
If $A,B,C$ are objects of $\Ac$ such that $F(A), F(B)\in\Bc$
and if $G_{AB}^C\neq\emptyset$, then $F(C)\in \Bc$ and $F$ is an 
isomorphism of complex varieties $G_{AB}^C \to G_{F(A), F(B)}^{F(C)}$. }
\qed

\vskip3mm

\noindent {\bf 3.2.} 
The characteristic function of 
a closed subvariety $W$ of an algebraic variety $Z$
is denoted by ${\bf 1}_W$. By ${\rm Fun}(Z)$ we denote the 
space of all constructible functions on $Z$. 
If $f\in{\rm Fun}(Z)$ its integral is: 
$$\int_Z f d\chi = \sum c_i \chi (W_i), \quad f=\sum c_i {\bf 1}_{W_i}.$$
Here $\chi$ is the Euler characteristic with compact support. 
More generally, if $\varphi: Z_1\to Z_2$ is any regular map
and $f\in {\rm Fun}(Z_1)$, then its direct image 
$\varphi_*(f)\in{\rm Fun}(Z_2)$ is (cf.[2]): 
$$(\varphi_*f)(z_2) = \int_{\phi^{-1}(z_2)} f d\chi.$$

\noindent
Similarly let $\Sc$ be an algebraic stacks of locally finite type.
The set of \Cb-points $\Sc(\Cb)$ 
is locally represented as the quotient of an algebraic variety by  
an action of an algebraic group.  
A constructible function on a stack $\Sc$ is a function $\Sc(\Cb)\to\Cb$ 
which can be 
represented as a finite linear combination of functions of 
the form  ${\bf 1}_{\Wc(\Cb)}$,
where $\Wc$ is a closed substack of finite type in $\Sc$. 
Let $\varphi: \Sc\to\Tc$ be a morphism of stacks  whose every fiber 
(over any \Cb-point) is an algebraic variety. 
Then we have the direct image map
$\varphi_*: {\rm Fun}(\Sc)\to{\rm Fun}(\Tc)$ defined as above. 

\vskip3mm

\noindent {\bf 3.3.} 
We now assume that $\Ac$ comes from a stack of Abelian categories 
over the category of
algebraic varieties, and that the moduli stack of objects of $\Ac$,
denoted by $\Ac^{iso}$,
is an algebraic stack of locally finite type. 

\vskip3mm

\noindent {\bf Examples.} 
(1) $\Ac$ is the category of representations of a finite-dimensional
{\bf C}-algebra. 

\noindent
(2) $\Ac$ is the category of coherent sheaves
on an algebraic variety with support in a projective subvariety.

\vskip3mm

\noindent
The set $\Ac^{iso}(\Cb)$ is the set of isomorphism classes of objects of 
$\Ac$.
The space ${\rm Fun}(\Ac^{iso})$ is made into an associative algebra, 
called the Hall algebra of $\cal A$ and
denoted by $\Hb(\Ac)$, as follows. 
Let $\Gc_\Ac$ be the stack formed by pairs $(A,B)$ of objects of $\Ac$,
with $A$ a subobject of $B$, and morphisms of such pairs.
There are three morphisms $p_1, p_2, p_3: \Gc_\Ac\to\Ac^{iso}$ 
which associate to $(A,B)$ the objects $A$, $B$, and $B/A$ respectively.
The fibers of $p_2$ are algebraic varieties.
The multiplication on $\Hb(\Ac)$ is 
$$f*g = p_{2*}((p_1^*f)\cdot(p_3^*g)).$$
Let $[A]\in\Hb(\Ac)$ be the characteristic function of the object $A$.
Then $\chi(G_{AB}^C)$ is the multiplicity of $[C]$ in $[A]*[B]$.

\vskip3mm

\noindent {\bf 3.4.}
Let $\Gamma$ be any finite graph without loops and multiple edges.
A double representation of $\Gamma$ is a rule which assigns
to each vertex $i$ a vector space $V_i$, and to any edge $\{i,j\}$ 
two operators $x_{ij}\,:\,V_j\to V_i$ and $x_{ji}\,:\,V_i\to V_j$ such that 
for every vertex $i$ we have $\sum_jx_{ij}x_{ji}=0.$
Finite dimensional double representations of $\Gamma$
form an Abelian category of finite type and global dimension 2 if $\Gamma$ is of affine type, 
denoted $\Rc_\Gamma$. 

\proclaim Proposition. 
Let $\Gamma$ be an affine Dynkin graph of type A-D-E, corresponding
to a finite subgroup $G\i SL_2(\Cb)$. 
The category $\Rc_\Gamma$ is equivalent to $\Cohc_G(\tau)$.

\vskip3mm

\noindent{\it Proof:} Let $\{\pi_i\}_{i\in I}$ be the set of simple
representations of $G$ and let $J=\{(i,j)\,|\,\pi_j\subset\pi_i\otimes\tau\}$
be the set of edges of the graph $\Gamma$. A $G$-equivariant coherent sheaf
on $\tau$ is a pair $(V,\phi)$ where $V$ is a finite dimensional $G$-module
and $\phi$ is a $G$-invariant linear map $V\otimes\Cb[\tau]\to V$.
Such a pair $(V,\phi)$ may be viewed as

$(i)$ a $I$-graded vector space $\bigoplus_{i\in I}V_{\pi_i}$,

$(ii)$ a collection of maps $\bigoplus_{(i,j)\in J}\,:\, 
V_{\pi_i}\otimes\tau\to\bigoplus_{(i,j)\in J}V_{\pi_j}.$

\noindent Then the relation $\sum_jx_{ij}x_{ji}=0$ means precisely that the 
maps
in $(ii)$ glue together in a map $V\otimes\Cb[\tau]\to V$.
\qed

\vskip3mm

\noindent 
Let $\Cb (i)$ be the simple object in $\Rc_\Gamma$
located at the vertex $i$ of $\Gamma$. 
Put $\theta_i=[\Cb(i)]\in\Hb(\Rc_\Gamma)$. 
As usual, set $\theta^{(k)}_i=\theta_i^k/k!$ for any $k\in\Nb^\times$.
Denote by $A_\Gamma$ the Cartan matrix of $\Gamma$, such that 
the index set is the set of vertices of $\Gamma$, $a_{ii}=2$, and $-a_{ij}$
is the number of edges joining $i$ and $j$ if $i\neq j$. 
Let ${\bf g}_\Gamma$ be the Kac-Moody Lie algebra associated to $A_\Gamma$. 
The ``nilpotent" subalgebra $\gb^+_\Gamma\i\gb_\Gamma$
is generated by the Chevalley generators 
$\eb_i$ subject to the Serre relations 
$$\forall i\neq j,\quad
\sum_k(-1)^k\eb_i^{\pm\,(k)}\eb_j^\pm\eb_i^{\pm\,(1-a_{ij}-k)}=0.$$
The following is a reformulation of a result from [8].

\vskip3mm

\noindent{\bf Theorem.} 
{\it The correspondence $\eb_i\to\theta_i$ defines a homomorphism 
$\Ub(\gb^+_\Gamma)\to \Hb(\Rc_\Gamma)$.} 
\qed

\vskip5mm

\centerline {\bf \S 4. Applications to algebraic surfaces.} 

\vskip3mm

\noindent
Let $S$ be a smooth surface over $\Cb$
and $C\i S$ be a reducible curve of which every component 
is a rational reduced $(-2)$-curve, such that
these curves meet transversely and not more than in one point. 
Let $\Pb^1_i, i\in I$, be the irreducible components of $C$.
Let $A_C$ be the negative of the intersection matrix of the components of $C$
and let $\gb_C$ be the Kac-Moody Lie algebra with Cartan matrix $A_C$.
Denote by $\Cohc(S,C)$ the category of coherent sheaves on $S$ 
with support in $C$. This category has finite type. 
Let $\Hb(S,C)$ be its Hall algebra. 

\vskip3mm

\noindent{\bf Theorem.} 
{\it The correspondence $\eb_i\mapsto [\Oc_{\Pb^1_i}]$
defines an algebra homomorphism $\Ub(\gb_C^+)\to\Hb(S,C)$.}

\vskip3mm
  
\proclaim Corollary. 
If $S$ is a projective surface then $\Ub(\gb_C^+)$ acts on 
${\rm Fun}(\Cohc(S))$. 
\qed

\vskip3mm

\noindent {\it Proof:} 
It is enough to consider the case when $C$ has type $A_2$. 
It is well-known (see [1], Theorem 7.3) that
if $C$ is a configuration of $(-2)$-curves on $S$ such 
that $A_C$ has type A-D-E
then the formal neighborhood of $C$ in $S$ is isomorphic to the formal 
neighborhood of the exceptional fiber $E$ in the Hilbert quotient 
$\Cb^2//G$ where $G\i SL_2(\Cb)$
is the finite subgroup with Cartan matrix $A_C$. 
Then we have the equivalence of categories 
$\Cohc(S,C) \simeq\Cohc(\Cb^2//G, E).$
Now, Theorem 1.4 and Propositions 3.1 and 3.4 imply 
that the subalgebra in $\Hb(\Rc_\Gamma)$ generated by the $\theta_\pi$,
$\pi\neq\Cb$, is isomorphic to the subalgebra in 
$\Hb(\Cb^2//G,E)$ generated by the $[\Oc_{\Pb^1_\pi}]$.
Thus  our theorem  follows from Theorem 3.4.
\qed

\vfill\eject 
\centerline {\bf References}

\vskip 10mm

\itemitem{[1]} W. Barth, C. Peters, A. Van de Ven, Compact Complex Surfaces,
{\it Springer} (1984). 

\itemitem{[2]} W. Fulton, R. MacPherson, Categorical Framework  for the
Study of Singular Spaces, {\it Memoirs AMS} {\bf 243}, 1981. 

\itemitem{[3]} G. Gonzalez-Sprinberg, J.L. Verdier, 
Construction g\'eom\'etrique de la correspondance
de McKay, {\it Ann. ENS}, {\bf 16} (1983), 409-449. 

\itemitem{[5]} Y. Ito, I. Nakamura, 
Hilbert schemes and simple singularities,
{\it New trends in algebraic geometry 
(Proceedings of the algebraic symposium, Warwick 1996)},
Cambridge University Press (1999), 151-233. 

\itemitem{[6]} M. Kapranov, Chow quotients of Grassmannians,
in: ``I.M. Gelfand Seminar" (S. Gelfand, S. Gindikin, Eds.), vol.1,
  ({\it Adv. in Soviet Math.}, v. 16, pt.2), p. 29-110, Amer. Math. Soc.
Providence RI 1993. 

\itemitem{[7]} P.B. Kronheimer, H. Nakajima, Yang-Mills instantons on
ALE gravitational instantons, {\it Math. Ann.} {\bf 288} (1990), 263-307

\itemitem{[8]} G. Lusztig, 
Quivers, perverse sheaves and enveloping algebras, 
{\it J. Amer. Math. Soc.}, {\bf 4} (1991), 365-421. 

\itemitem{[9]} J. McKay, Graphs, singularities and finite groups, in: 
``The Santa Cruz conference on finite groups' (Santa Cruz CA 1979), p.
183-186, Pro. Symp. Pure Math. {\bf 37}, Amer. Math. Soc. Providence RI
1980.

\itemitem{[10]} H. Nakajima, Instantons on ALE spaces, quiver varieties
and Kac-Moody algebras, {\it Duke Math. J. } {\bf 76} (1994), 365-416. 

\itemitem{[11]} H. Nakajima,
Lectures on Hilbert schemes of points on surfaces,
{\it Preprint} (1998).

\itemitem{[12]} C.M. Ringel, Hall algebras and quantum groups, {\it
Invent. Math.} {\bf 101} (1990), 583-592.

\vskip1cm

\noindent
M.K. : Department of Mathematics, Northwestern University,\hfill\break
Evanston, Illinois 60208.
kapranov at math.nwu.edu.

\vskip5mm

\noindent
E.V. : D\'epartement de Math\'ematiques, Universit\'e de Cergy-Pontoise,
\hfill\break
2 Av. A. Chauvin, 95302 Cergy-Pontoise Cedex, France.
vasserot at math.pst.u-cergy.fr.

\bye